\documentclass{IEEEtran4PSCC}
\pdfoutput=1

\newcommand{\ba}{\begin{array}}
\newcommand{\ea}{\end{array}}
\newcommand{\bc}{\begin{center}}
\newcommand{\ec}{\end{center}}

\newcommand{\beqn}[1]{\begin{equation}\label{#1}}
\newcommand{\eeqn}{\end{equation}}
\newcommand{\be}{\begin{equation}}
\newcommand{\ee}{\end{equation}}

\newcommand{\beqnn}{\begin{eqnarray}}
\newcommand{\eeqnn}{\end{eqnarray}}

\newcommand{\bs}{\boldsymbol}
\newcommand{\tp}{\mathsf{T}}

\usepackage{amsmath}
\usepackage{graphicx}
\usepackage{amssymb}
\usepackage{cite}
\usepackage{cases}

\usepackage{multicol}
\usepackage[figuresright]{rotfloat}[2004/01/04]
\usepackage{amsfonts}
\usepackage{setspace}
\usepackage{cite}
\usepackage{mathrsfs}
\usepackage{graphics,dblfloatfix} 
\usepackage{epsfig} 
\usepackage{mathdots}
\usepackage[usenames, dvipsnames]{color}
\allowdisplaybreaks[1]
\usepackage{algorithm}
\usepackage{algpseudocode}
\usepackage{soul}
\usepackage{xargs}                      
\usepackage[pdftex,dvipsnames]{xcolor}  
\usepackage[colorinlistoftodos,prependcaption,textsize=scriptsize, textwidth=10mm]{todonotes}
\newcommandx{\revfirst}[2][1=]{\todo[linecolor=blue,backgroundcolor=blue!25,bordercolor=blue,#1]{#2}}
\newcommandx{\revsecond}[2][1=]{\todo[linecolor=red,backgroundcolor=red!25,bordercolor=red,#1]{#2}}
\newcommandx{\revthird}[2][1=]{\todo[linecolor=Orange,backgroundcolor=Orange!25,bordercolor=Orange,#1]{#2}}
\newcommandx{\thiswillnotshow}[2][1=]{\todo[disable,#1]{#2}}

\IEEEoverridecommandlockouts

\begin{document}
\title{Chance-Constrained Shrunken-Primal-Dual Subgradient (CC-SPDS) Approach for Decentralized Electric Vehicle Charging Control}

\author{
\IEEEauthorblockN{Mingxi Liu, Mostafa Sahraei-Ardakani}
\IEEEauthorblockA{Department of Electrical and Computer Engineering \\
University of Utah\\
Salt Lake City, UT, USA 84112\\
Email: mingxi.liu@utah.edu}
\thanks{\textcopyright 2019 IEEE.  Personal use of this material is permitted.  Permission from IEEE must be obtained for all other uses, in any current or future media, including reprinting/republishing this material for advertising or promotional purposes, creating new collective works, for resale or redistribution to servers or lists, or reuse of any copyrighted component of this work in other works.}
}


\maketitle

\begin{abstract}
In this paper, we develop a chance-constrained decentralized electric vehicle (EV) charging control framework to achieve ``valley-filling'' meanwhile meeting individual charging requirements and satisfying distribution network constraints. The control design is formulated as an optimization problem with a stochastic non-separable objective function and globally coupled chance constraints. We propose a novel chance-constrained shrunken-primal-dual subgradient (CC-SPDS) algorithm to support the chance-constrained decentralized control scheme and verify its efficacy and convergence with a representative distribution network model.
\end{abstract}



\section{Introduction} \label{Introduction}

Electric vehicles (EVs) have shown their effectiveness in increasing energy conversion efficiency and reducing GHG emissions \cite{Liu_PIEEE_2013}. Besides, with proper control, the charging process of a large fleet of EVs can be leveraged for provisioning grid-level services \cite{Callaway_PIEEE_2011}. Extensive studies have been conducted on the modeling and control problems of EVs for ``valley-filling'' \cite{Gan_TPS_2013}, load balancing \cite{Mercurio_MCCA_2013}, and frequency regulation \cite{Karfopoulos_TPS_2016_2}. The potential of EV charging control for facilitating renewable energy integration and providing additional services is discussed in \cite{Richardson_RSER_2013}. 

Apart from the benefits, charging a large fleet of EVs without proper control would have negative impacts on distribution networks. Lopes \emph{et al.} \cite{Lopes_PIEEE_2011} surveyed the challenges of integrating EVs into the mid- and low-voltage distribution networks, including voltage sag, transformer overloading, network congestion, and increased power loss. The negative impacts caused by non-network-aware controlled charging were verified in \cite{Fernandez_TPS_2011}. Hence, an impact-free EV integration necessitates designing EV charging control in the context of distribution network constraints. 

Recently, an increasing amount of attention is being paid to controlling EV charging for distribution network impact alleviation and grid service provisioning. Considering network impacts, Richardson \emph{et al.} \cite{Richardson_TPS_2012} optimized EV charging profiles to minimize the total power consumption. In \cite{Luo_IET_2013}, a real-time control for power loss minimization was studied. Authors in \cite{Zhang_TPS_2017} and \cite{Liu_TCST_2017} separately developed decentralized EV charging controllers to achieve valley-filling under voltage constraints. Studies in a large amount of literature have shown that one of the most effective services EVs could provide is to fill the load valley at night \cite{Ma_TCST_2013}. In our paper, we target at valley-filling meanwhile meeting both local charging needs and distribution network constraints on voltage.


Most existing EV charging control algorithms assume a powerful central controller that can handle large computational loads, \cite{Clement_Nyns_TPS_2010, Richardson_TPS_2012, Luo_IET_2013, Quiros-Tortos_TPS_2016}. As the EV deployment in the distribution network keeps increasing, centralized approaches however do not scale well and become unrealistic. In this paper, we aim to design an optimization-based decentralized control scheme that does not require a communication network among EV chargers. Achieving this is notably challenging, because firstly the valley-filling objective constructs a coupled and non-separable function, and secondly the network voltage constraints strongly couple the individual charging powers in form of linear inequalities \cite{Liu_TCST_2017}. For this type of problem, Koshal \emph{et al.} \cite {Koshal_SIAM_2011} developed a regularized primal-dual subgradient approach, however convergence errors exist due to the regularization of the Lagrangian. In \cite{Zhang_TPS_2017}, the authors proposed an Alternating Direction Method of Multipliers (ADMM)-based decentralized algorithm. Though proved to be practically applicable, the required two-layer communication network complicates the communication and poses computing burdens to all buses. Additionally, ADMM-based algorithms often face the problem of a large number of iterations. In our previous work \cite{Liu_TCST_2017}, we proposed a shrunken-primal-dual subgradient (SPDS) approach that can eliminate convergence errors, reduce the number of iterations, and alleviate communication loads. 

Unfortunately, the aforementioned algorithms are not capable of handling uncertainties in EV charging (e.g., random baseline load and random customer's behaviors). Among the limited literature, Liu \emph{et al.} \cite{Liu_CDC_2017} proposed an SPDS-based event-triggering decentralized control algorithm to handle EV drivers' random arrivals and departures, where customers' baseline loads are assumed perfectly predicted. Hassan \emph{et al.} \cite{Hassan_PSCC_2018} developed a chance-constrained ADMM for decentralized distributed energy resources control considering random solar generation. However, decentralized implementation of chance constraints on voltage was treated deterministically, leaving that problem open. Handling chance constraints on voltage in a decentralized fashion is notably challenging because: (i) Voltage constraints are strongly coupled; (ii) Gradient of the chance constraint function is ill-conditioned, leading to zero updates in iteration-based algorithms. In this paper, we aim to leverage the outstanding convergency, optimality, and speed of SPDS, to design a chance-constrained SPDS (CC-SPDS) to handle decentralized EV charging control problems that have chance constraints on nodal voltages. To the author's best knowledge, this has never been attempted before.

The main contribution of this paper is two-fold. First, this is the first paper studying decentralized EV charging control for provisioning valley-filling under chance constraints due to uncertain baseline loads. Second, a novel CC-SPDS algorithm is developed to solve a class of optimization problems consisting of non-separable objective functions and strongly coupled network chance constraints.

\section{Preliminaries \& Problem formulation} \label{Problem_formulation}

\subsection{EV charging model} \label{general_charging_model}
For a radial distribution network with $h$ nodes, let $n_\imath$ denote the number of EVs connected at the $\imath$th node and $n=\sum_{\imath=1}^{h}n_\imath$ denote the total number of EVs. The individual EV charging dynamics can be represented as \cite{Liu_CDC_2017}
\begin{equation} \label{E_left_model}
x_{\imath,\hat{\imath}}(T+1)= x_{\imath,\hat{\imath}}(T)+B_{\imath,\hat{\imath}}u_{\imath,\hat{\imath}}(T), \nonumber
\end{equation}
where the subscript $_{\imath,\hat{\imath}}$ denotes the $\hat{\imath}$th EV connected at the $\imath$th node,  $T$ is the discrete-time index, $x_{\imath,\hat{\imath}}(T)$ denotes the energy remaining to be charged, $B_{\imath,\hat{\imath}}{=}-\eta_{\imath,\hat{\imath}}\Delta t \bar{P}_{\imath,\hat{\imath}}$, $\eta_{\imath,\hat{\imath}}$ is the charging efficiency, $\Delta t$ is the sampling interval, $\bar{P}_{\imath,\hat{\imath}}$ is the maximum charging power, and $0 {\leq} u_{\imath,\hat{\imath}}(T) {\leq} 1$ is the control signal.

Let $k_{\imath,\hat{\imath}}$ and $k_{\imath,\hat{\imath}}{+}K_{\imath,\hat{\imath}}$ denote individual plug-in time and designated charging deadlines; let $[k, k{+}K]$ denote the valley-filling service period. In this paper, we assume that individual charging periods $[k_{\imath,\hat{\imath}}, k_{\imath,\hat{\imath}}{+}K_{\imath,\hat{\imath}}]$ are contracted to fully cover the valley-filling period. This assumption can be easily relaxed. The charging efficiencies, initial state of charge, maximum charging powers, and battery capacities are heterogeneous. Augmenting all connected EVs, we have the system dynamics represented as
\begin{equation} \label{state_space_state_all}
x(T+1)=x(T)+\sum_{\imath=1}^{h}\sum_{\hat{\imath}=1}^{n_\imath}B_{{\imath,\hat{\imath}}}^cu_{\imath,\hat{\imath}}(T),
\end{equation}
where $x(T)=[x_{1,1}(T)\cdots x_{1,n_1}(T) \cdots x_{h,1}(T)\cdots x_{h,n_h}(T)]^{\mathsf{T}}$, $u(T)=[u_{1,1}(T)\cdots u_{1,n_1}(T) \cdots u_{h,1}(T)\cdots u_{h,n_h}(T)]^{\mathsf{T}}$, and $B_{\imath,\hat{\imath}}^c$ is the $(\imath,\hat{\imath})$th column of the diagonal matrix $\text{diag}\left\{ B_{\imath,\hat{\imath}}\right\} \in \mathbb{R}^{n \times n}$, $\hat{\imath}=1,\ldots,n_\imath,~\imath=1,\ldots,h$.

Augmenting the system state $x(T)$ and individual control signal $u_{\imath,\hat{\imath}}(T)$ in \eqref{state_space_state_all} along the valley-filling period $[k,k+K]$, we have 
\begin{equation}
\begin{aligned}
\mathcal{X}(k)&=\left[ x(k+1|k)^{\mathsf{T}}~x(k+2|k)^{\mathsf{T}}~\cdots~x(k+K|k)^{\mathsf{T}} \right]^{\mathsf{T}}, \\
\mathcal{U}_{\imath,\hat{\imath}}(k)&=\left[ u_{\imath,\hat{\imath}}(k|k)~u_{\imath,\hat{\imath}}(k+1|k) \cdots u_{\imath,\hat{\imath}}(k+K-1|k)\right]^{\mathsf{T}}. \nonumber
\end{aligned}
\end{equation}
Herein, $x(k+\kappa|k)$, $\kappa=1,\cdots,K$, is the system state at time $k+\kappa$ predicted at time $k$; $u_{\imath,\hat{\imath}}(k+\kappa-1|k)$, $\kappa=1,\cdots,K$, is the control signal at time $k+\kappa-1$ predicted at time $k$.


To meet all drivers' charging requirements, the energy remaining to be charged at the end of valley-filling, i.e., the $K$th vector block in $\mathcal{X}(k)$, must satisfy
\begin{equation} \label{Charging_requirement_constraint}
x(k+K|k)=x(k)+\sum_{\imath=1}^h\sum_{\hat{\imath}=1}^{n_\imath}\mathcal{B}_{{\imath,\hat{\imath}}}^l\mathcal{U}_{\imath,\hat{\imath}}(k) =\boldsymbol{0},
\end{equation}
where $\mathcal{B}_{{\imath,\hat{\imath}}}^l=[ B_{{\imath,\hat{\imath}}}^c~B_{{\imath,\hat{\imath}}}^c~\cdots~B_{{\imath,\hat{\imath}}}^c] \in \mathbb{R}^{n \times K}$.

\subsection{Distribution network model} \label{network_model_section}
In this paper, we consider radial distribution networks. Let $\mathbb{H}=\{\imath| \imath=1,\ldots,h\}$ denote the set of downstream nodes and let $\mathbb{S}$ denote the set of all downstream line segments. The feeder head decouples interactions in the downstream distribution system from the rest of the grid and maintains its own voltage magnitude $\left| V_0 \right|$. To better illustrate the EV charging impacts and how controlled EV charging can help alleviate the impacts, we do not consider any distributed energy resources (DERs), voltage regulations, or reactive power supplies.

At time $T$, let $\left| V_\imath (T) \right|$ denote the voltage magnitude at Node $\imath$; let $p_\imath(T)$ and $q_\imath(T)$ denote the real and reactive power consumption at Node $\imath$; and with a slight abuse of notations, let $r_{\imath \jmath}+ix_{\imath \jmath}$ denote the impedance of the line segment $(\imath, \jmath)$. By omitting line losses, the LinDistFlow model of this distribution network is derived as \cite{Baran_TPD_1989, Farivar_CDC_2013, Liu_TCST_2017}
\begin{equation} \label{LinDistFlow_original}
\boldsymbol{V}(T)=\boldsymbol{V}_0-2\boldsymbol{R}p(T)-2\boldsymbol{X}q(T),
\end{equation}
where
\begin{equation}
\begin{aligned}
\boldsymbol{V}(T)&=[ \left| V_1(T) \right|^2~\left| V_2(T) \right|^2~\cdots~\left| V_h(T) \right|^2 ]^{\mathsf{T}}, \\
\boldsymbol{V}_0&= [ \left| V_0 \right|^2 ~ \left| V_0 \right|^2~\cdots~\left| V_0 \right|^2 ]^{\mathsf{T}}, \\
p(T)&=\left[ p_1(T)~p_2(T)~\cdots~p_h(T)\right]^{\mathsf{T}} , \\
q(T)&=\left[ q_1(T)~q_2(T)~\cdots~q_h(T)\right]^{\mathsf{T}} , \nonumber
\end{aligned}
\end{equation}
and
\begin{equation} \label{Definition_R_X}
\begin{aligned}
&\boldsymbol{R}\in\mathbb{R}^{h\times h},~\boldsymbol{R}_{\imath \jmath}=\sum_{(\hat{\imath},\hat{\jmath})\in \mathbb{S}_\imath \cap \mathbb{S}_\jmath}r_{\hat{\imath}\hat{\jmath}}, \\
&\boldsymbol{X}\in\mathbb{R}^{h\times h},~\boldsymbol{X}_{\imath \jmath}=\sum_{(\hat{\imath},\hat{\jmath})\in \mathbb{S}_\imath \cap \mathbb{S}_\jmath}x_{\hat{\imath}\hat{\jmath}}, \nonumber
\end{aligned}
\end{equation}
where $\mathbb{S}_\imath$ ($\mathbb{S}_\jmath$) is the set containing downstream line segments connecting the feeder head and Node $\imath$ (Node $\jmath$) \cite{Farivar_CDC_2013}.

At each node and each house, real and reactive loads can be separated into baseline loads and EV charging loads, indicating 
\begin{equation}
\begin{aligned}
p_\imath(T)&=\sum_{\imath=1}^{n_\imath}\left(p_{\imath, \hat{\imath}}^b(T)+p_{\imath,\hat{\imath}}^{EV}(T) \right), \\
q_\imath(T)&=\sum_{\imath=1}^{n_\imath}\left(q_{\imath, \hat{\imath}}^b(T)+q_{\imath,\hat{\imath}}^{EV}(T) \right),
\end{aligned} \nonumber
\end{equation}
where $p_{\imath, \hat{\imath}}^b(T)$ and $q_{\imath, \hat{\imath}}^b(T)$ denote the baseline real and reactive power consumption, respectively, and $p_{\imath,\hat{\imath}}^{EV}(T)$ and $q_{\imath,\hat{\imath}}^{EV}(T)$ denote the real and reactive EV charging power, respectively. As it is unrealistic to assume perfect forecast of the baseline load, a realistic problem formulation must consider $p_{\imath, \hat{\imath}}^b(T){\sim} \mathcal{N} (\mu^{p}_{\imath, \hat{\imath}}(T), \sigma^{p}_{\imath, \hat{\imath}}(T))$ and $q_{\imath, \hat{\imath}}^b(T){\sim} \mathcal{N} (\mu^{q}_{\imath, \hat{\imath}}(T), \sigma^{q}_{\imath, \hat{\imath}}(T))$, where $\mu^{p}_{\imath, \hat{\imath}}(T)$ and $\mu^{q}_{\imath, \hat{\imath}}(T)$ denote the means of the normal distributions, and $\sigma^{p}_{\imath, \hat{\imath}}(T)$ and $\sigma^{q}_{\imath, \hat{\imath}}(T)$ denote the standard deviations of the normal distributions.
In this paper, we make the following three assumptions:
\begin{itemize}
\item[{\bf{A1}}] At each time $T$, power consumption of all houses follows independent normal distributions. The standard derivation $\sigma^{p}_{\imath, \hat{\imath}}(T)$ is assumed to be a constant $\sigma^{p}$ at all time. 
\item[{\bf{A2}}] EVs do not consume reactive power, i.e., $q_{\imath,\hat{\imath}}^{EV}(T)=0$.
\item[{\bf{A3}}] Power factor at each house is a constant, which implies $q_{\imath, \hat{\imath}}^b(T) = \gamma p_{\imath, \hat{\imath}}^b(T)$. 
\end{itemize}
\noindent{\bf{Remark 1:}} In {\bf{A1}}, other distributions can also be used, but normal distribution can better reflect customers' random behaviors. Constant standard deviation is assumed for the simplicity of presentation. {\bf{A2}} follows as no fourth-quadrant charger is considered. {\bf{A3}} follows by assuming power factor correction
equipment is available at each node. \hfill $\blacksquare$

Let $\hat{p}^b(T)=[p_{1,1}^b(T)\cdots p_{h,n_h}^b(T) ]^\tp \in \mathbb{R}^n$ contain baseline power at all $n$ houses, then $\hat{p}^b(T) {\sim} \mathcal{N}\left( \bs{\mu}^p(T), \bs{\Sigma}^p\right)$, where $ \bs{\mu}^p(T) = [ \mu^{p}_{1, 1}(T) \cdots \mu^{p}_{h, n_h}(T)]^\tp$ and $\bs{\Sigma}^p=\text{diag} \left\{ (\sigma^p)^2 \right\} \in \mathbb{R}^{n \times n}$. Further let $p^b(T)=G \hat{p}^b(T)=[p_1^b(T)~\cdots~p_h^b(T)]^\tp$ denote the nodal aggregated baseline power, where $G=\text{diag}\left\{ G_\imath \right\}\in \mathbb{R}^{h \times n}$, $\imath=1,\ldots, h$, $G_\imath={\bold{1}}_{n_\imath}^\mathsf{T}$ is the nodal aggregation vector. Since \eqref{LinDistFlow_original} is linear, the voltage drop caused by the baseline load can be represented as
\begin{equation}
\boldsymbol{V}_b(T)=2\bs{R}Gp^b(T)+2\bs{X}Gq^b(T) = \bs{H}p^b(T), \nonumber
\end{equation}
where $\bs{H}=2(\bs{R}+\gamma \bs{X})G$, yielding
\begin{equation}
\boldsymbol{V}(T)=\boldsymbol{V}_0-\boldsymbol{V}_b(T)-2\boldsymbol{R}G\bar{P}u(T), \nonumber
\end{equation}
where $\bar{P}=\text{diag}\left\{ \bar{P}_{\imath,\hat{\imath}} \right\}$, $\hat{\imath}=1,\ldots,n_\imath$, $\imath=1,\ldots,h$ and $\boldsymbol{V}_b(T) \sim \mathcal{N}(  \tilde{\bs{\mu}}^p(T),  \tilde{\bs{\Sigma}}^p )$ with $ \tilde{\bs{\mu}}^p(T) = \bs{H} \bs{\mu}^p(T)$ and $\tilde{\bs{\Sigma}}^p =\bs{H} \bs{\Sigma}^p \bs{H}^\tp$.


Let $D\in \mathbb{R}^{h\times n}$ denote $-2\boldsymbol{R}G\bar{P}$, $y_d(T)$ denote $\boldsymbol{V}_0-\boldsymbol{V}_b(T)$, and $y(T)$ denote $\boldsymbol{V}(T)$, we have
\begin{equation}
y(T)=y_d(T)+Du(T), \nonumber
\end{equation}
where $y_d(T) {\sim} \mathcal{N} ( \bs{V}_0-\tilde{\bs{\mu}}^p(T), \tilde{\bs{\Sigma}}^p )$. The system output augmented along the valley-filling period is written as
\begin{equation}
\mathcal{Y}_{k}=\mathcal{Y}_{dk}+\sum_{\imath=1}^{h}\sum_{\hat{\imath}=1}^{n_\imath}\mathcal{D}_{\imath,\hat{\imath}}\mathcal{U}_{\imath,\hat{\imath}}(k), \nonumber
\end{equation}
where
\begin{equation}
\begin{aligned}
\mathcal{Y}_{k}&=\left[y(k|k)^{\mathsf{T}}~ y(k+1|k)^{\mathsf{T}}~\cdots y(k+K-1|k)^{\mathsf{T}} \right]^{\mathsf{T}}, \\
\mathcal{Y}_{dk}&=\left[y_d(k|k)^{\mathsf{T}}~ y_d(k+1|k)^{\mathsf{T}}~\cdots y_d(k+K-1|k)^{\mathsf{T}} \right]^{\mathsf{T}}, \\
\mathcal{D}_{\imath,\hat{\imath}}&=\bigoplus_{\kappa=1}^{K}D_{\imath,\hat{\imath}}, ~D=\left[ D_{1,1}\cdots D_{1,n_1} \cdots D_{h,1} \cdots D_{h,n_h}\right] \nonumber
\end{aligned}
\end{equation}
and
\begin{equation}
\mathcal{Y}_{dk} {\sim} \mathcal{N} \left( \left[ \begin{array}{c}
\bs{V}_0-\tilde{\bs{\mu}}^p(k) \\
\vdots \\
\bs{V}_0-\tilde{\bs{\mu}}^p(k+K-1)
\end{array} \right],  \left[ \begin{array}{ccc}
 \tilde{\bs{\Sigma}}^p & ~ & ~ \\
~ & \ddots & ~ \\
~ & ~ & \tilde{\bs{\Sigma}}^p
\end{array}\right] \right). \nonumber
\end{equation}
Herein, $\oplus$ denotes the direct matrix sum.

\subsection{Deterministic valley-filling} \label{deterministic_valley_filling}
Suppose $p^b_{\imath, \hat{\imath}}(T)$ is deterministic or can be accurately forecasted. According to \cite{Liu_TCST_2017, Liu_CDC_2017}, optimal EV charging sequences for valley-filling under nodal voltage constraints are obtained by solving
\begin{equation} \label{original_problem}
\begin{aligned}
\min_{\mathcal{U}} &~\mathcal{F}(\mathcal{U})=\frac{1}{2}\left\| P_b+\tilde{P}\mathcal{U} \right\|_2^2+\frac{\rho}{2}\left\| \mathcal{U}\right\|_2^2 \\
\text{s.t.} &~\mathcal{U}_{\imath,\hat{\imath}}\in \mathbb{U}_{\imath,\hat{\imath}},~\forall~\imath=1,\ldots,h,~\hat{\imath}=1,\ldots,n_{\imath}, \\
&~d(\mathcal{U})\triangleq\underline{\nu}^2\hat{\boldsymbol{V}}_0-\mathcal{Y}_{dk}-\sum_{\imath=1}^{h}\sum_{\hat{\imath}=1}^{n_\imath}\mathcal{D}_{\imath,\hat{\imath}}\mathcal{U}_{\imath,\hat{\imath}} \leq {\bold{0}},
\end{aligned}
\end{equation}
where the time index $k$ is dropped for simplicity, $\hat{\boldsymbol{V}}_0=\left[{\boldsymbol{V}}_0^\tp \cdots {\boldsymbol{V}}_0^\tp \right]^\tp \in \mathbb{R}^{nK}$ and
\begin{equation}
\begin{aligned}
\mathcal{U}&=[\mathcal{U}_{1,1}^{\mathsf{T}}\ldots \mathcal{U}_{1,n_1}^{\mathsf{T}} \ldots \mathcal{U}_{h,1}^{\mathsf{T}}\ldots \mathcal{U}_{h,n_h}^{\mathsf{T}}]^{\mathsf{T}}, \\
\mathbb{U}_{\imath,\hat{\imath}} &:= \left\{ \mathcal{U}_{\imath,\hat{\imath}} | {\bold{0}}\leq \mathcal{U}_{\imath,\hat{\imath}} \leq {\bold{1}}, x_{\imath,\hat{\imath}}(k)+\mathcal{B}_{\imath, \hat{\imath}}^l\mathcal{U}_i=0\right\}, \nonumber
\end{aligned}
\end{equation}
Herein, $P_b=\left[ \left\| p^b(k)  \right\|_1 \cdots \left\| p^b(k+K-1)  \right\|_1 \right]^\tp$ denotes the aggregated baseline load profile at the feeder head, $\tilde{P}$ is the EV charging power aggregation matrix, $\rho/2 \left\| \mathcal{U}\right\|_2^2$ is for battery state of health protection, and $\underline{\nu}$ is the lower bound of the transformer service range.

\subsection{Chance-constrained valley-filling} \label{stochastic_valley_filling}
In more realistic situations, where baseline load cannot be accurately forecasted, we have normally distributed $p_{\imath, \hat{\imath}}^b(T){\sim} \mathcal{N} (\mu^{p}_{\imath, \hat{\imath}}(T), \sigma^{p})$. Under this circumstance, the deterministic problem \eqref{original_problem} becomes
\begin{equation} \label{stochastic_problem}
\begin{aligned}
\min_{\mathcal{U}} &~\underset{p^b_{\imath,\hat{\imath}}{\sim}\mathcal{N}(\cdot,\cdot)}{\mathbb{E}}\left[ \mathcal{F}(\mathcal{U}) \right] \\
\text{s.t.} &~\mathcal{U}_{\imath,\hat{\imath}}\in \mathbb{U}_{\imath,\hat{\imath}},~\forall~\imath=1,\ldots,h,~\hat{\imath}=1,\ldots,n_{\imath}, \\
&~\text{Pr}\left( [ z ]_\kappa \geq [\hat{\mathcal{Y}} ]_\kappa  \right) \geq \delta,~\forall~\kappa=1,\ldots, K,\\ 
\end{aligned}
\end{equation}
where $[\hat{\mathcal{Y}} ]_\kappa {\sim} \mathcal{N} ( (\underline{\nu}^2-1)\bs{V}_0+\tilde{\bs{\mu}}^p(k+\kappa-1), \tilde{\bs{\Sigma}}^p ) \in \mathbb{R}^h$ is the $\kappa$th block of $\hat{\mathcal{Y}}=\underline{\nu}^2 \hat{\bs{V}}_0-\mathcal{Y}_{dk}$, $[ z ]_\kappa= [ \sum_{\imath=1}^{h}\sum_{\hat{\imath}=1}^{n_\imath}\mathcal{D}_{\imath,\hat{\imath}}\mathcal{U}_{\imath,\hat{\imath}}]_\kappa \in \mathbb{R}^{h}$ is the $\kappa$th block of $ \sum_{\imath=1}^{h}\sum_{\hat{\imath}=1}^{n_\imath}\mathcal{D}_{\imath,\hat{\imath}}\mathcal{U}_{\imath,\hat{\imath}}$, and $1-\delta$ is the allowed chance of violating voltage limits. For simplicity of presentation, we use $\mathbb{E}$ to denote the expectation operator $\mathbb{E}_{p^b_{\imath,\hat{\imath}}{\sim}\mathcal{N}(\cdot,\cdot)}$ hereinafter. Designing decentralized algorithms for stochastic problem \eqref{stochastic_problem} is more complicated than for problem \eqref{original_problem} as the coupled objective function becomes an expectation and the coupled inequality constraints become coupled chance constraints.

\section{Chance-constrained SPDS design} \label{stochastic_SPDS_design}
Let the multivariate normal cumulative distribution function (MVNCDF) of $[\hat{\mathcal{Y}}]_\kappa$ be denoted by $F_{[\hat{\mathcal{Y}}]_\kappa}(\cdot)$, then the chance constraints in problem \eqref{stochastic_problem} become
\begin{equation} \label{chance_d}
d_\kappa(\mathcal{U})=-F_{[\hat{\mathcal{Y}}]_\kappa}([z]_\kappa)+\delta \leq 0,~\forall~\kappa=1,\ldots, K,
\end{equation}
yielding the unconstrained Lagrangian of problem \eqref{stochastic_problem} as
\begin{equation}
\mathcal{L}(\mathcal{U},\bs{\lambda})={\mathbb{E}}\left[ \mathcal{F}(\mathcal{U}) \right] + \sum_{\kappa=1}^{K}\lambda_\kappa d_\kappa(\mathcal{U}). \nonumber
\end{equation}

It is proved in \cite{Liu_TCST_2017} that a convex, strongly coupled optimization problem can be efficiently solved by SPDS by following
\begin{small}
\begin{subequations} \label{constant_update}
\begin{equation} \label{primal_update}
\mathcal{U}_{\imath,\hat{\imath}}^{(\ell+1)} {=}\Pi_{\mathbb{U}_{\imath,\hat{\imath}}}\left(\frac{1}{\tau_\mathcal{U}}\Pi_{\mathbb{U}_{\imath,\hat{\imath}}}\left( \tau_\mathcal{U}\mathcal{U}_{\imath,\hat{\imath}}^{(\ell)}{-}\alpha\nabla_{\mathcal{U}_{\imath,\hat{\imath}}}\mathcal{L}(\mathcal{U}^{(\ell)},\bs{\lambda}^{(\ell)})\right)\right),
\end{equation}
\begin{equation} \label{dual_update}
{\lambda_{\kappa}}^{(\ell+1)}{=}\Pi_{\mathbb{D}}\left(\frac{1}{\tau_{\lambda_{\kappa}}}\Pi_{\mathbb{D}}\left(\tau_{\lambda_\kappa} {\lambda_{\kappa}}^{(\ell)}{+}\beta\nabla_{{\lambda_{\kappa}}}\mathcal{L}(\mathcal{U}^{(\ell)},\bs{\lambda}^{(\ell)})\right)\right),~~
\end{equation}
\end{subequations}
\end{small}where $\Pi_{\cdot}(\cdot)$ is the Euclidean projection, $\ell$ denotes the iteration number, $\alpha$ and $\beta$ are the update step sizes for primal and dual variables, respectively, $0<\tau_{\mathcal{U}}, \tau_{\lambda_\kappa}<1$ are the shrinking parameters, and $\mathbb{D}$ is the constraint set for dual variables. The definition of $\mathbb{D}$ is referred to \cite{Liu_TCST_2017, Liu_CDC_2017}. It has been shown in \cite{Liu_TCST_2017} that the convergence and optimality are guaranteed if $\alpha$ and $\beta$ are appropriately chosen for a designated tuple $(\tau_{\mathcal{U}}, \tau_{\lambda_\kappa})$.

Nevertheless, SPDS in \eqref{constant_update} cannot be directly applied to the chance-constrained problem \eqref{stochastic_problem} due to two challenges: (i) The feasible region is non-convex due to the non-convexity of $F_{[\hat{\mathcal{Y}}]_\kappa}(\cdot)$; (ii) The gradient $\nabla_{\mathcal{U}_{\imath,\hat{\imath}}}d_\kappa(\mathcal{U})$ is bounded and in most cases it realizes the value of $0$, implying that updates of the primal variables can barely capture constraints. (This will be illustrated later). In the following, we will present the details of CC-SPDS and show via simulations that it can handle the non-convex chance-constrained problem \eqref{stochastic_problem}. Rigorous proof of convergency and optimality gap will be provided in a journal extension of this paper.

The objective function of problem \eqref{stochastic_problem} can be rewritten as
\begin{equation}
\begin{aligned}
{\mathbb{E}}\left[ \mathcal{F}(\mathcal{U}) \right] &= \frac{1}{2}{\mathbb{E}} \left[  \left\| P_b+\tilde{P}\mathcal{U} \right\|_2^2 \right] +\frac{\rho}{2}\left\| \mathcal{U}\right\|_2^2 \\
&\leftrightarrow  {\mathbb{E}}^\tp \left[  P_b  \right] (\tilde{P}\mathcal{U}) + \frac{1}{2}\left\| \tilde{P}\mathcal{U}\right\|_2^2+\frac{\rho}{2}\left\| \mathcal{U}\right\|_2^2,
\end{aligned} \nonumber
\end{equation}
and its gradient at $\mathcal{U}^{(\ell)}$ can be calculated as
\begin{equation}
\begin{aligned}
\nabla_{\mathcal{U}_{\imath,\hat{\imath}}} \mathbb{E}[\mathcal{F}(\mathcal{U}^{(\ell)})] =& \bar{P}_{\imath,\hat{\imath}}\left(\left[ \begin{array}{c}
\left\| \bs{\mu}^p(k)\right\|_1 \\
\vdots \\
\left\| \bs{\mu}^p(k+K-1)\right\|_1 
\end{array} \right] + \tilde{P}\mathcal{U}^{(\ell)} \right) \\
&+\rho \mathcal{U}_{\imath,\hat{\imath}}^{(\ell)},
\end{aligned} \nonumber
\end{equation}
where inside the bracket of the first term on the right hand side is the mean of the total load profile at the $\ell$th iteration. This is a universal information broadcasted to all EVs in each iteration. To obtain the complete form of $\nabla_{\mathcal{U}_{\imath,\hat{\imath}}}\mathcal{L}(\mathcal{\mathcal{U},\bs{\lambda}})$ we first introduce the following theorem.

\noindent{\bf{Theorem 1}}\cite{Prekopa_sto_book} {\em Let $\xi {\sim} \mathcal{N}(\bs{\mu},\bs{\Sigma})$ with some positive definite covariance matrix $\bs{\Sigma}=(\sigma_{ij}) \in \mathbb{R}^{s \times s}$. Then, the distribution function $F_\xi$ is continuously differentiable at any $z\in \mathbb{R}^s$ and
\begin{equation}
\frac{\partial F_\xi}{\partial z_j}(z)=f_{\xi_j}(z_j)\cdot F_{\xi_{\bar{j}}}(z_{\{1,\ldots,s \}\setminus j}). \nonumber
\end{equation}
Here, $f_{\xi_j}$ denotes the 1-D normal density of the component $\xi_j$, $\xi_{\bar{j}}$ is an $(s-1)$-D normal random vector distributed according to $\xi_{\bar{j}}{\sim} \mathcal{N}(\hat{\bs{\mu}}, \hat{\bs{\Sigma}})$, $\hat{\bs{\mu}}$ results from the vector $\bs{\mu}+\sigma_{jj}^{-1}(z_j-\mu_j)\bs{\sigma}_j$ by deleting component $j$, and $\hat{\bs{\Sigma}}$ results from the matrix $\bs{\Sigma}-\sigma_{jj}^{-1}\bs{\sigma}_j\bs{\sigma}_j^\tp$ by deleting row $j$ and column $j$, where $\bs{\sigma}_j$ refers to column $j$ of $\bs{\Sigma}$. Moreover, $\hat{\bs{\Sigma}}$ is positive definite. \hfill $\blacksquare$}

By applying Theorem 1, we can write
\begin{equation} \label{original_d_k}
\begin{aligned}
\nabla_{\mathcal{U}_{\imath,\hat{\imath}}}d_\kappa(\mathcal{U})&= -\nabla_{\mathcal{U}_{\imath,\hat{\imath}}} F_{[\hat{\mathcal{Y}}]_\kappa}([z]_\kappa) \\
&=-\bs{J}^\tp_{\mathcal{U}_{\imath,\hat{\imath}}}([z]_\kappa) \nabla_{[z]_k}F_{[\hat{\mathcal{Y}}]_\kappa}([z]_\kappa) \\
&=-\mathcal{D}_{\imath,\hat{\imath}}^\tp G_\kappa^\tp   \left[ \begin{array}{c}
 f_{[\hat{\mathcal{Y}}]_{\kappa,1}}([z]_{\kappa,1}) \cdot F_{[\hat{\mathcal{Y}}]_{\kappa,\bar{1}}}([z]_{\kappa,\bar{1}}) \\
\vdots \\
 f_{[\hat{\mathcal{Y}}]_{\kappa,h}}([z]_{\kappa,h}) \cdot F_{[\hat{\mathcal{Y}}]_{\kappa,\bar{h}}}([z]_{\kappa,\bar{h}}) 
\end{array}\right],
\end{aligned}
\end{equation}
where $G_\kappa = [ \bs{0}\cdots \bs{0}~ \underbrace{\bs{I}_{h}}_{\kappa\text{th}}\ ~\bs{0}\cdots \bs{0} ]$, $\bs{J}(\cdot)$ denotes the Jacobian matrix, $[z]_{\kappa, j}$ denotes the $j$th component of $[z]_{\kappa}$, and $[z]_{\kappa, \bar{j}}$ is obtained by removing the $j$th component of $[z]_{\kappa}$.\

Though theoretically sound, directly applying $\nabla_{\mathcal{U}_{\imath,\hat{\imath}}}d_\kappa(\mathcal{U})$ in \eqref{original_d_k} to the primal update equation \eqref{primal_update} is ill-conditioned. For example in the independent bivariate normal distribution case, the gradient of the cumulative function is $\nabla F_{x,y}(X,Y) = [f_x(X)F_y(Y)~f_y(Y)F_x(X)]^\tp$. \begin{figure}[!htb] \centering
\includegraphics[width=0.42\textwidth, trim = 25mm 2mm 25mm 3mm, clip]{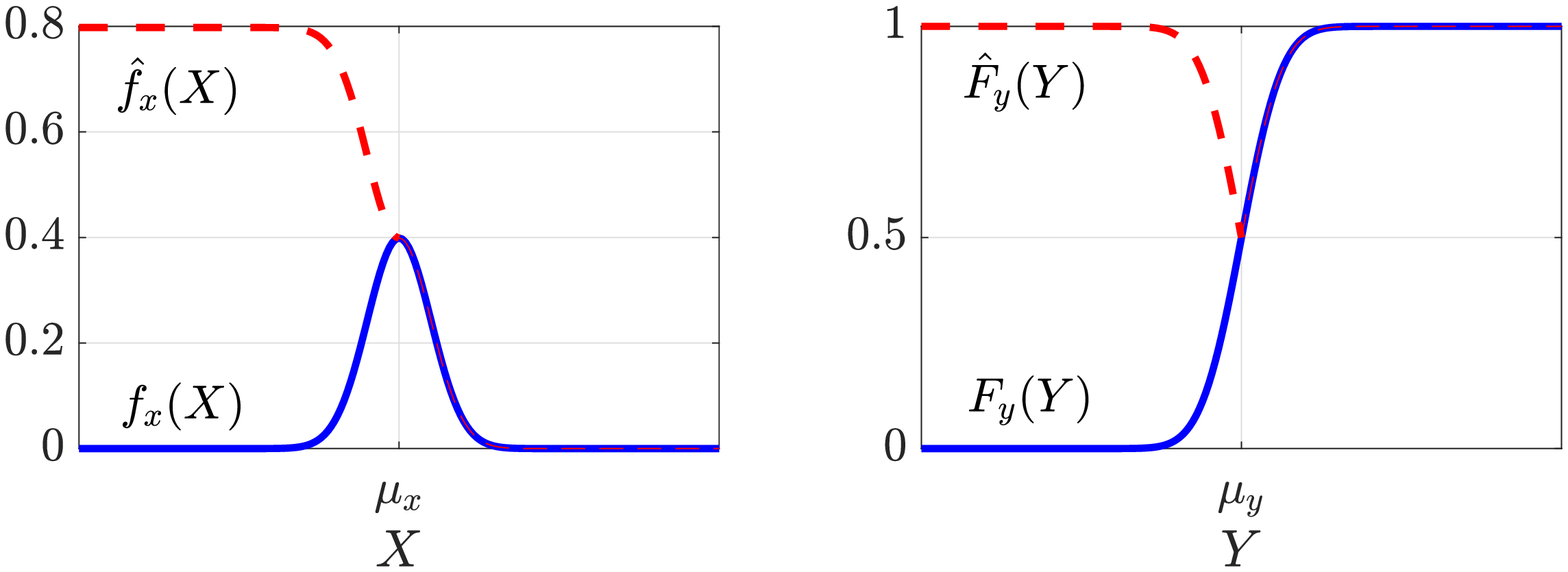}
\caption{Illustration of the ill-conditioned gradient of MVNCDF.}
\label{pdf_cdf_illus}
\end{figure} As illustrated in Figure \ref{pdf_cdf_illus}, there are four cases where $f_x(X)F_y(Y)$ is extremely small or zero: (i \& ii) $X$ is far away from $\mu_x$ to the right and $Y$ is far away from $\mu_y$ to the left (and right); (iii \& iv) $X$ is far away from $\mu_x$ to the left and $Y$ is far away from $\mu_y$ to the left (and right); Translated into \eqref{chance_d}, these four cases mean: (i \& ii) Constraint is satisfied with a high probability; No (or slight) updates w.r.t. the constraint are needed in the iteration; (iii \& iv) Constraint is not satisfied with a high probability; Updates w.r.t. the constraint are needed in the iteration. However, in the last two cases, updates cannot be realized as the updating direction (gradient) is $0$. To overcome this ill-conditioned gradient, we define an adjusted gradient as
\begin{equation}
\hat{\nabla} F_{x,y}(X,Y) \triangleq [\hat{f}_x(X)\hat{F}_y(Y)~\hat{f}_y(Y)\hat{F}_x(X)]^\tp, \nonumber
\end{equation}
where
\begin{equation}
\begin{aligned}
& \hat{f}_x(X)\triangleq\left\{ \begin{array}{ll}
f_x(X) &~X\geq \mu_x \\
2f_x(\mu_x)-f_x(2\mu_x-X) &~X< \mu_x 
\end{array}\right. \\
&\hat{F}_y(Y)\triangleq\left\{ \begin{array}{ll}
F_y(Y) &~~~~~~~~~~~~~~Y\geq \mu_y \\
F_y(2\mu_y-Y) &~~~~~~~~~~~~~~Y< \mu_y 
\end{array}\right.
\end{aligned} \nonumber
\end{equation}
and $\hat{f}_y(Y)$, $\hat{F}_x(X)$ are defined in the same way. The adjusted density function and cumulative function are illustrated in Figure \ref{pdf_cdf_illus} as the dashed curves which only adjust the left half planes. These adjustments guarantee effective updates once chance constraints are violated. Following this idea, we have ${\nabla}_{\mathcal{U}_{\imath,\hat{\imath}}}d_\kappa(\mathcal{U})$ in \eqref{original_d_k} adjusted  to
\begin{equation} \label{adjusted_d}
\hat{\nabla}_{\mathcal{U}_{\imath,\hat{\imath}}}d_\kappa(\mathcal{U})=-\mathcal{D}_{\imath,\hat{\imath}}^\tp G_\kappa^\tp   \left[ \begin{array}{c}
 \hat{f}_{[\hat{\mathcal{Y}}]_{\kappa,1}}([z]_{\kappa,1}) \cdot \hat{F}_{[\hat{\mathcal{Y}}]_{\kappa,\bar{1}}}([z]_{\kappa,\bar{1}}) \\
\vdots \\
 \hat{f}_{[\hat{\mathcal{Y}}]_{\kappa,h}}([z]_{\kappa,h}) \cdot \hat{F}_{[\hat{\mathcal{Y}}]_{\kappa,\bar{h}}}([z]_{\kappa,\bar{h}}) 
\end{array}\right],
\end{equation}
and consequently the primal update equation \eqref{primal_update} becomes
\begin{equation} \label{primal_update_adjusted}
\mathcal{U}_{\imath,\hat{\imath}}^{(\ell+1)} =\Pi_{\mathbb{U}_{\imath,\hat{\imath}}}\left(\frac{1}{\tau_\mathcal{U}}\Pi_{\mathbb{U}_{\imath,\hat{\imath}}}\left( \tau_\mathcal{U}\mathcal{U}_{\imath,\hat{\imath}}^{(\ell)}-\alpha\hat{\nabla}_{\mathcal{U}_{\imath,\hat{\imath}}}\mathcal{L}(\mathcal{U}^{(\ell)},\bs{\lambda}^{(\ell)})\right)\right),
\end{equation}
where
\begin{equation}
\hat{\nabla}_{\mathcal{U}_{\imath,\hat{\imath}}}\mathcal{L}(\mathcal{U}^{(\ell)},\bs{\lambda}^{(\ell)})=\nabla_{\mathcal{U}_{\imath,\hat{\imath}}} \mathbb{E}[\mathcal{F}(\mathcal{U}^{(\ell)})]+\sum_{\kappa=1}^{K}\lambda_\kappa \hat{\nabla}_{\mathcal{U}_{\imath,\hat{\imath}}}d_\kappa(\mathcal{U}^{(\ell)}). \nonumber
\end{equation}
The dual update equation \eqref{dual_update} can be readily calculated as
\begin{equation}
\nabla_{{\lambda_{\kappa}}}\mathcal{L}(\mathcal{U}^{(\ell)},\bs{\lambda}^{(\ell)})=d_\kappa(\mathcal{U}^{(\ell)}). \nonumber
\end{equation}
Together \eqref{primal_update_adjusted} and \eqref{dual_update} construct the CC-SPDS.

\noindent{\bf{Remark 2:}} In EV charging control, the adjusted cumulative gradient vector on the right hand side of \eqref{adjusted_d} is broadcasted to all EVs as a \emph{public key}. Each EV uses its \emph{private key} $\mathcal{D}_{\imath,\hat{\imath}}$ to interpret the public key, so that the control sequence update can fully acknowledge the global chance constraints. \hfill $\blacksquare$

\noindent{\bf{Remark 3:}} CC-SDPS, especially the approach of adjusting the gradient of MVNCDFs, can be applied to many other popular distributions, e.g., Poisson distribution and uniform distribution, with $\hat{f}$ and $\hat{F}$ properly defined. \hfill $\blacksquare$


\section{Simulations} \label{Simulations}
The simplified IEEE 13 Node Test Feeder used in \cite{Liu_TCST_2017} is adopted in this paper. In specific, Node 632 and Node 671 are not connected to any load, and each of the other nodes is connected with $70$ houses equipped with level-2 chargers. Battery capacities are uniformly distributed in $[18, 20]$ kWh. Initial and designated SOCs are uniformly distributed in $[0.3,0.5]$ and $[0.7,0.9]$, respectively. Primal step size $\alpha=1\times10^{-11}$; dual step size $\beta=2$; $d_{\lambda_{\kappa}}=5\times 10^5$ in $\mathbb{D}$; shrinking parameters are $\tau_\mathcal{U}=\tau_{\lambda_\kappa}=0.974$. The above parameters were empirically chosen to accelerate the convergence speed. Initial values of $\mathcal{U}_i^{(0)}$ and $\lambda_{\kappa}^{(0)}$ are set to zero vectors. The voltage lower bound is chosen to $\underline{\nu}=0.954$, which is slightly higher than the $0.95$ p.u. in ANSI C84.1 service standard, to compensate for the discarded line losses in LinDistFlow model and possible voltage violations. For the normal distribution of baseline load at each house, the mean is obtained by scaling the data from Southern California Edison \cite{SCE_reg}; the constant standard deviation $\sigma^p=400$W. The contracted valley-filling service period is from 19:00 to 8:00 next day. The sampling time is 15 minutes. The allowed voltage constraint violation chance $1{-}\delta$ is set to $10\%$.

\begin{figure}[!htb] \centering
\includegraphics[width=0.35\textwidth, trim = 15mm 2mm 23mm 4mm, clip]{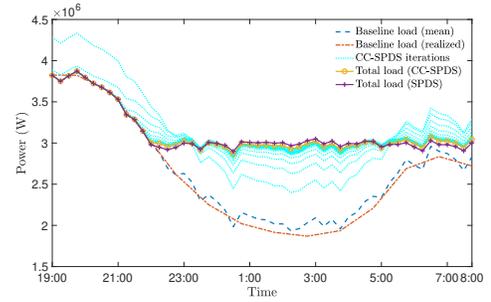}
\caption{Valley-filling performance of CC-SPDS and SPDS.}
\label{valley_filling_performance}
\end{figure}

Figure \ref{valley_filling_performance} shows the valley-filling performance under the proposed CC-SPDS decentralized framework, where convergence of the algorithm can be readily revealed. In total 50 iterations were executed and each iteration took about 0.16 second. Comparing with the purely deterministic case discussed in \cite{Liu_TCST_2017, Liu_PSCC_2018} where the converged total load profile is flat, the converged total load profile in the stochastic case presents slight variations. This is because the objective function only captures the mean of the baseline load profile, i.e., CC-SPDS fills the mean baseline load valley, and once random behaviors realized the valley-filling performance is slightly compromised. Figure \ref{valley_filling_performance} also shows the total load profile resulting from directly applying deterministic SPDS to the stochastic case, i.e., solving \eqref{original_problem} by using the mean baseline loads in both the objective function and voltage constraints. The difference between CC-SPDS and SPDS w.r.t. valley-filling performance is neglectable, as both of them are using mean baseline load in the objective function.

\begin{figure}[!htb] \centering
\includegraphics[width=0.35\textwidth, trim = 0mm 0mm 0mm 12mm, clip]{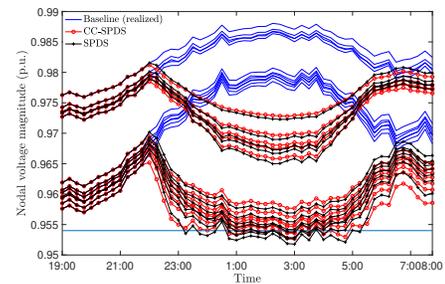}
\caption{Nodal voltage magnitudes from baseline load, CC-SPDS controlled total load, and SPDS controlled total load.}
\label{voltage_magnitude}
\end{figure}

Figure \ref{voltage_magnitude} shows the nodal voltage magnitudes of all nodes during the valley-filling period. Under baseline loads, all nodal voltage magnitudes are above $0.954$ p.u. lower bound. Under the proposed CC-SPDS-based decentralized EV charging control framework, only few nodes occasionally exceed the bound during 1:00 to 5:00. Under the SPDS-based control, voltage magnitudes are continuously violating the limit during 23:00 to 5:30 and are much lower than those of the CC-SPDS control. This observation reveals the effectiveness of CC-SPDS. 

\begin{figure}[!htb] \centering
\includegraphics[width=0.38\textwidth, trim = 0mm 16mm 0mm 13mm, clip]{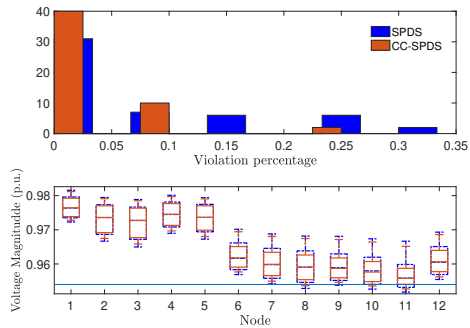}
\caption{Upper: Histogram of violation events at every 15 minutes; Lower: Box plots of nodal voltage magnitudes during valley-filling, where solid are CC-SPDS and dot dashed are SPDS.}
\label{voltage_violation_comp}
\end{figure}

Figure \ref{voltage_violation_comp} presents the statistics of voltage violations.  The upper plot reveals that in the total 52 time slots during the valley-filling period, CC-SPDS can suppress 50 (96\%) of them to have less than 10\% (only one) of all nodes violating 0.954 p.u., and 2 (4\%) of them to have 25\% (three) of all nodes violating 0.954 p.u. At each time slot, percentage of the total nodes violating the limit can be treated as a measure of the chance constraint. The 25\% violation is as expected since the chance constraints only reflect a probability which cannot be fully reflected by a single sample. Comparing with CC-SPDS, SPDS has the histogram more distributed to the right, leading to more violation events and greater violation magnitude. In the lower plot, boxes represent the 25th and 75th percentiles and legs with bars represent the extreme values. It can be concluded that under the control of CC-SPDS: Nodal voltage magnitudes have less variations; The 25th percentiles are guaranteed to stay above 0.954 p.u.; Number and magnitudes of violations are less than those of SPDS. These observations further reveal the effectiveness of CC-SPDS.

\section{Conclusion} \label{Conclusion}
This paper developed a novel decentralized EV charging control framework where distribution network constraints were formulated as chance constraints due to customers' random behaviors. A new decentralized optimization algorithm -- CC-SPDS was developed to solve the chance-constrained valley-filling problem as well as generic optimization problems that have coupled objective functions and coupled global chance constraints. Convergence and effectiveness of CC-SPDS were verified via simulations. CC-SPDS can be generalized to other distributions and extended to consider local chance constraints.

\bibliographystyle{IEEEtran}
%

\end{document}